\documentclass{amsart}
\usepackage{amsmath,amssymb,amsthm,mathrsfs}

\newcommand{\BC}{{\mathbb C}}\newcommand{\BD}{{\mathbb D}}

\newcommand{\cD}{{\mathcal D}}
\newcommand{\cE}{{\mathcal E}}\newcommand{\cF}{{\mathcal F}}
\newcommand{\cG}{{\mathcal G}}\newcommand{\cH}{{\mathcal H}}

\newcommand{\cK}{{\mathcal K}}\newcommand{\cL}{{\mathcal L}}
\newcommand{\cM}{{\mathcal M}}\newcommand{\cN}{{\mathcal N}}

\newcommand{\cR}{{\mathcal R}}

\newcommand{\cU}{{\mathcal U}}\newcommand{\cV}{{\mathcal V}}
\newcommand{\cW}{{\mathcal W}}\newcommand{\cX}{{\mathcal X}}
\newcommand{\cY}{{\mathcal Y}}

\newcommand{\bH}{{\mathbf H}}

\newcommand{\bS}{{\mathbf S}}

\newcommand{\sL}{{\mathscr L}}

\newcommand{\sR}{{\mathscr R}}

\newcommand{\sV}{{\mathscr V}}
\newcommand{\sW}{{\mathscr W}}

\newcommand{\tilV}{\tilde{V}}

\newcommand{\wtilA}{\widetilde{A}}\newcommand{\wtilB}{\widetilde{B}}
\newcommand{\wtilC}{\widetilde{C}}\newcommand{\wtilD}{\widetilde{D}}

\newcommand{\wtilH}{\widetilde{H}}

\newcommand{\wtilK}{\widetilde{K}}

\newcommand{\wtilQ}{\widetilde{Q}}\newcommand{\wtilR}{\widetilde{R}}
\newcommand{\wtilT}{\widetilde{T}}
\newcommand{\wtilU}{\widetilde{U}}\newcommand{\wtilV}{\widetilde{V}}
\newcommand{\wtilW}{\widetilde{W}}

\newcommand{\Ga}{\Gamma}
\newcommand{\De}{\Delta}

\newcommand{\Th}{\Theta}

\newcommand{\la}{\lambda}\newcommand{\La}{\Lambda}

\newcommand{\vph}{\varphi}
\newcommand{\om}{\omega}

\newcommand{\im}{\textup{Im\,}}

\newcommand{\kr}{\textup{Ker\,}}

\newcommand{\mat}[2]{\ensuremath{\left[\begin{array}{#1}#2\end{array} \right]}}

\newcommand{\sbm}[1]{\left[\begin{smallmatrix} #1\end{smallmatrix}\right]}
\newcommand{\ov}[1]{{\overline{#1}}}

\newcommand{\tu}[1]{\textup{#1}}

\newcommand{\Llra}{\Longleftrightarrow}

\newcommand{\half}{\frac{1}{2}}
\newcommand{\ands}{\quad\mbox{and}\quad}

\newcommand{\ons}{\mbox{ on }}
\newcommand{\LDS}{\{A,T',U',R,Q\}}
\newcommand{\tilom}{\widetilde\omega}
\newcommand{\cFt}{\widetilde\cF}

\newcommand{\SO}{X}                      
\newcommand{\wtilSO}{\widetilde\SO}

\newcommand{\CS}{\underline{0}}

\newtheorem{theorem}{Theorem}[section]
\newtheorem{corollary}[theorem]{Corollary}
\newtheorem{lemma}[theorem]{Lemma}
\newtheorem{proposition}[theorem]{Proposition}

\begin{document}

\title{Redheffer representations and relaxed commutant lifting}

\author{S. ter Horst}
\address{Department of Mathematics\\
Utrecht University\\
Budapestlaan 6\\
3584 CD Utrecht\\
The Netherlands}
\email{terhorst@vt.edu}

\begin{abstract}
It is well known that the solutions of a (relaxed) commutant lifting
problem can be described via a linear fractional representation of
the Redheffer type. The coefficients of such Redheffer
representations are analytic operator-valued functions defined on
the unit disc $\BD$ of the complex plane. In this paper we consider
the converse question. Given a Redheffer representation, necessary
and sufficient conditions on the coefficients are obtained
guaranteeing the representation to appear in the description of the
solutions to some relaxed commutant lifting problem. In addition, a
result concerning a form of non-uniqueness appearing in the
Redheffer representations under consideration and an harmonic
maximal principle, generalizing a result of A. Biswas, are proved.
The latter two results can be stated both on the relaxed commutant
lifting as well as on the Redheffer representation level.
\end{abstract}



\maketitle

\setcounter{section}{-1}\setcounter{equation}{0}
\section{Introduction}\label{S:intro}

Linear fractional representations have been used for describing solutions
to metric constrained interpolation problems since the work of
Adamyan-Arov-Kre\u\i n \cite{AAK68,AAK71}, and later appeared in the encompassing commutant
lifting theory of Sz.-Nagy-Foia\c s \cite{NF68} and D. Sarason \cite{S67};
cf., \cite{FF90,FFGK98}. They also play an important role in linear system theory \cite{ZDG96}.
The topic of the present paper is the class of linear fractional representation that appears
in the context of the description of the solutions to the relaxed commutant lifting problem
obtained in \cite{tH07}. Such representations are called Redheffer representations, after the
work of  R. Redheffer \cite{R60,R62}.

Before we can introduce this class of representations, we require some definitions and notations.
Let $\cU$ and $\cY$ be Hilbert spaces. With the symbol $\bH^\infty(\cU,\cY)$
we denote the set of uniformly bounded analytic functions on the open
unit disc $\BD$ with values in $\sL(\cU,\cY)$. Here, as usual,
we write $\sL(\cU,\cY)$ for the space of operators from $\cU$ into $\cY$. By definition, an
operator is assumed to be linear and bounded.
The set $\bH^\infty(\cU,\cY)$ is a Banach space with respect to the supremum
norm $\|\ \|_\infty$; its closed unit ball is
denoted by $\bS(\cU,\cY)$ and called the {\em Schur class} associated with $\cU$
and $\cY$. Functions in $\bS(\cU,\cY)$ are called {\em Schur class functions}.
It is well known that a function $F$ in $\bH^\infty(\cU,\cY)$ defines a multiplication
operator $M_F$ between the Hardy spaces $H^2(\cU)$ and $H^2(\cY)$, and that $\|F\|_\infty$ is equal to the
operator norm of $M_F$ in $\sL(H^2(\cU),H^2(\cY))$.
We further define $\bH^2(\cU,\cY)$ to be the set of analytic functions
$H$ on $\BD$ whose values are in $\sL(\cU,\cY)$ with the property that the formula
\begin{equation}\label{GaH}
(\Ga_H u)(\la)=H(\la)u\quad(u\in\cU,\la\in\BD)
\end{equation}
defines an operator from $\cU$ into the Hardy space $H^2(\cY)$. The set $\bH^2(\cU,\cY)$
is a Banach space with norm $\|H\|=\|\Ga_{H}\|$, i.e., the
norm of $H\in\bH^2(\cU,\cY)$ is equal to the operator norm of the associated operator
$\Ga_H$. We write $\bH^2_\tu{ball}(\cU,\cY)$ for the closed unit ball
of $\bH^2(\cU,\cY)$. Note that, conversely, any operator $\Ga\in\sL(\cU,H^2(\cY))$
defines a function $H\in\bH^2(\cU,\cY)$ via the same identity $H(\la)u=(\Ga u)(\la)$.

Now let $\cU$, $\cY$, $\cE$ and $\cE'$ be Hilbert spaces and $\Psi$ an operator-valued
function on $\BD$ that  decomposes as
$\Psi=\sbm{\Psi_{1,1}&\Psi_{1,2}\\ \Psi_{2,1}&\Psi_{2,2}}$ so that $\Psi_{1,1}(0)=0$ and
\begin{equation}\label{RedCoefs}
\Psi_{1,1}\in\bS(\cE',\cE),\ \
\Psi_{1,2}\in\bH^2_\tu{ball}(\cU,\cE),\ \
\Psi_{2,1}\in\bS(\cE',\cY),\ \
\Psi_{2,2}\in\bH^2_\tu{ball}(\cU,\cY),
\end{equation}
The {\em Redheffer representation} associated with $\Psi$ is the map
$V\mapsto\sR_\Psi[V]$ given by
\begin{equation}\label{RedTrans}
\sR_\Psi[V](\la)=
\Psi_{2,2}(\la)+\Psi_{2,1}(\la)V(\la)
(I-\Psi_{1,1}(\la)V(\la))^{-1}\Psi_{1,2}(\la)
\end{equation}
and defined for $V\in\bS(\cE,\cE')$. The functions $\Psi_{1,1}$, $\Psi_{1,2}$, $\Psi_{2,1}$ and
$\Psi_{2,2}$ are referred to as the associated {\em Redheffer coefficients}. In addition, we shall
assume the {\em coefficient matrix}
\begin{equation}\label{RedMat}
\mat{cc}{M_{\Psi_{1,1}}&\Ga_{\Psi_{1,2}}\\M_{\Psi_{2,1}}&\Ga_{\Psi_{2,2}}}
\mat{c}{H^2(\cE')\\\cU}\to\mat{c}{H^2(\cE)\\H^2(\cY)}
\end{equation}
to be a contraction, or, more often a co-isometry. It is then the case that the Redheffer representation
$\sR_\Psi$ maps $\bS(\cE,\cE')$ into $\bH^2_\tu{ball}(\cU,\cY)$ (see Proposition \ref{P:Redheffer}).
Note that this includes the more conventional class of Redheffer representations where the function
$\Psi$ is a Schur class function itself, in which case $\cR_\Psi$ maps $\bS(\cE,\cE')$ into $\bS(\cU,\cY)$.

The relaxed
commutant lifting problem was introduced in \cite{FFK02}; the developed theory extends
the classical commutant lifting theory \cite{NF68}, as well as the Treil-Volberg
commutant lifting setting \cite{TV94} and the weighted version of \cite{BFF99}.
See \cite{DBF} for an application in filterbank design.
In \cite{FFK02} a particular (central) solution is given; while descriptions of all solutions
are obtained in \cite{LT06,FtHK06a,FtHK06b,tH07,tH1}. The starting point for the relaxed
commutant lifting problem is a {\em lifting data set} $\LDS$ consisting
of five Hilbert space operators: the operator $A$ is a contraction mapping $\cH$ into ${\cH}'$, the
operator $U'$ on $\cK'$ is a minimal isometric lifting of the contraction $T'$ on $\cH'$, i.e.,
$U'$ is an isometry on $\cK'$ with $\cH'\subset\cK'$ being cyclic for $U'$ and
$\Pi_{\cH'}U'=T'\Pi_{\cH'}$, and $R$ and $Q$ are operators from $\cH_0$ to $\cH$, satisfying
\begin{equation}\label{intertw}
T'AR=AQ\ands R^*R\leq Q^*Q.
\end{equation}
Given this data set the {\em relaxed commutant lifting problem} is to describe a (all) contraction(s)
$B$ from $\cH$ to $\cK'$ such that
\begin{equation}\label{rclt}
\Pi_{\cH'}B=A\ands U'BR=BQ.
\end{equation}
Here we follow the convention that for a subspace $\cV$ of a Hilbert space $\cW$ the symbol
$\Pi_\cV$ stands for the orthogonal projection from $\cW$ onto $\cV$ viewed as an operator
from $\cW$ onto $\cV$, whereas $P_\cV$ is used for the orthogonal projection from $\cW$ onto $\cV$
as an operator from $\cW$ into $\cW$.
A contraction $B$ from $\cH$ into $\cK'$ that satisfies \eqref{rclt} is called a
{\em contractive interpolant for $\LDS$}. The Treil-Volberg version  appears when $R$ is the identity
operator on $\cH$, and thus $\cH_0=\cH$; for classical commutant lifting it is assumed in addition
that $Q$ is an isometry.

Without loss of generality we may, and will, assume that the isometric lifting $U'$ in
$\LDS$ is the Sz.-Nagy Sch\"affer isometric lifting of $T'$, i.e.,
\[
U'=\mat{cc}{T'&0\\E_{\cD_{T'}}D_{T'}&S_{\cD_{T'}}}\ons\mat{c}{\cH'\\H^2(\cD_{T'})},
\]
where we use the general notation $S_\cU$ for the forward shift on the Hardy space $H^2(\cU)$,
and $E_{\cU}$ for the canonical embedding of $\cU$ into the subspace of constant functions in
$H^2(\cU)$, that is, $(E_\cU u)(\la)=u$ for all $\la\in\BD=\{\la\in\BC\colon |\la|<1\}$ and $u\in\cU$.
Furthermore, as usual, given a contraction $N$, we write $D_N$ for the defect operator and $\cD_N$
for the defect space of $N$, that is, $D_N$ is the positive square root of $I-N^*N$ and
$\cD_N$ is the closure of the range of $D_N$.

An initial step in the process of describing the solutions to a commutant lifting problem is the
extraction of an operator from the data. In the case of the relaxed commutant lifting problem
we obtain a contraction $\om$ defined by
\begin{equation}\label{om}
\om:\cF=\ov{D_AQ\cH_0}\to\mat{c}{\cD_{T'}\\\cD_A},\quad \om
D_AQ=\mat{c}{D_{T'}AR\\D_AR\\},
\end{equation}
which we refer to as the {\em underlying contraction} of $\LDS$ (see \cite{FFK02}).
We write $\om_1$ for the component of $\om$ that maps $\cF$ into $\cD_{T'}$
and $\om_2$ for the component of $\om$ that maps $\cF$ into $\cD_{A}$.
In the classical commutant lifting setting, as well as in many of the examples considered
in \cite{FFK02}, we have $R^*R=Q^*Q$, which is equivalent to $\om$ being an isometry.

Once the underlying contraction $\om$ is obtained from the data,  we
set $\cG=\cD_A\ominus\cF$ and write $\Pi_\cF$ and $\Pi_\cG$ for the
orthogonal projections in $\sL(\cD_A,\cF)$ and $\sL(\cD_A,\cG)$,
respectively. The symbols $\Pi_{\cD_A}$ and $\Pi_{\cD_{T'}}$ will be
used for the orthogonal projections from $\cD_A\oplus\cD_{T'}$ onto
$\cD_A$, respectively $\cD_{T'}$. Next we define operator-valued
functions $\Phi_{1,1}$, $\Phi_{1,2}$, $\Phi_{2,1}$ and $\Phi_{2,2}$
on $\BD$ by
\begin{equation}\label{RCLcoefs}
\begin{array}{rl}
&\Phi_{1,1}(\la)=\la\Pi_\cG(I_{\cD_A}-\la\om_2\Pi_\cF)^{-1}
\Pi_{\cD_A} D_{\om^*},\\[.1cm]
&\Phi_{1,2}(\la)=\Pi_\cG(I_{\cD_A}-\la\om_2\Pi_\cF)^{-1},\\[.1cm]
&\Phi_{2,1}(\la)=\Pi_{\cD_{T'}} D_{\om^*}+\la\om_1\Pi_\cF
(I_{\cD_A}-\la\om_2\Pi_\cF)^{-1}\Pi_{\cD_A} D_{\om^*},\\[.1cm]
&\Phi_{2,2}(\la)=\om_1\Pi_\cF(I_{\cD_A}-\la\om_2\Pi_\cF)^{-1}.
\end{array}\quad(\la\in\BD)
\end{equation}
and put
\begin{equation}\label{PhiRCL}
\Phi=\mat{cc}{\Phi_{1,1},&\Phi_{1,2}\\\Phi_{2,1}&\Phi_{2,2}}.
\end{equation}
The set of all contractive interpolants is then described by the following
theorem; see Theorem 5.1.1 in \cite{tH07}.

\begin{theorem}\label{T:RCLRedheffer}
Let $\LDS$ be a lifting data set with underlying contraction $\om$.
Set $\cG=\cD_A\ominus D_AQ\cH_0$ and define $\Phi_{1,1}$, $\Phi_{1,2}$,
$\Phi_{2,1}$ and $\Phi_{2,2}$ by \eqref{RCLcoefs}. Then
\begin{equation*}
\Phi_{1,1}\!\in\bS(\cD_{\om^*},\cG),\
\Phi_{1,2}\!\in\bH^2_\tu{ball}(\cD_A,\cG),\
\Phi_{2,1}\!\in\bS(\cD_{\om^*},\cD_{T'}),\
\Phi_{2,2}\!\in\bH^2_\tu{ball}(\cD_A,\cD_{T'}),
\end{equation*}
$\Phi_{1,1}(0)=0$ and the coefficient matrix
\begin{equation}\label{CoefMat}
K_0=\mat{cc}{M_{\Phi_{1,1}}&\Ga_{\Phi_{1,2}}\\M_{\Phi_{2,1}}&\Ga_{\Phi_{2,2}}}
:\mat{c}{H^2(\cD_{\om^*})\\\cD_A}\to\mat{cc}{H^2(\cG)\\ H^2(\cD_{\om^*})}
\end{equation}
is a co-isometry. Moreover, set
$H_V=\sR_\Phi[V]\in\bH^2_\tu{ball}(\cD_A,\cD_{T'})$
for any $V$ in $\bS(\cG,\cD_{\om^*})$ and define
\begin{equation}\label{RCLsols2}
B_V=\mat{c}{A\\\Ga_{H_V} D_A}:\cH\to\mat{c}{\cH'\\ H^2(\cD_{T'})}.
\end{equation}
Then $B_V$ is a contractive interpolant for $\LDS$, and all contractive
interpolants are obtained in this way. Finally, $K_0$ is unitary if and
only if $\om$ is an isometry and $\om_2\Pi_\cF$ on $\cD_A$ is strongly
stable (i.e., $\lim_{n\to\infty}(\om_2\Pi_\cF)^n u=0$ for each $u\in\cD_A$).
\end{theorem}

When we take for $V$ the constant function whose value is the
zero operator, we see that a particular contractive interpolant is obtained
by taking $H_V=\Phi_{2,2}$ in \eqref{RCLsols2}; this contractive interpolant is
called the {\em central contractive interpolant} and was already obtained
in \cite{FFK02}.

One point where the Redheffer representation of Theorem
\ref{T:RCLRedheffer} is distinctively different from the one in the
classical commutant lifting setting is that the map $\sR_\Phi$
is, in general, not one-to-one, and thus the same holds true for the
map $V\mapsto B_V$. This property is inherited from the Schur
representation of \cite{FtHK06b} (from which the Redheffer
representation of Theorem \ref{T:RCLRedheffer} is deduced in
\cite{tH07}) where the same phenomenon occurs. We return to this
issue later in the introduction.

Theorem \ref{T:RCLRedheffer} raises the question how the set of Redheffer representations
$\sR_\Phi$ with $\Phi$ obtained from a relaxed commutant lifting problem is situated in the
set of all Redheffer representations $\sR_\Psi$ with coefficients as in \eqref{RedCoefs}.
We further restrict the set of Redheffer representations $\sR_\Psi$ by demanding that
$\Psi_{1,1}(0)=0$ and that the coefficient matrix \eqref{RedMat} is a co-isometry, as this is
a natural restriction based on the result of Theorem \ref{T:RCLRedheffer}.

To be precise, we will consider the following question:~Given a Redheffer representation
$\cR_\Psi$ with coefficients as a in \eqref{RedCoefs}, with $\Psi_{1,1}(0)=0$ and \eqref{RedMat}
a co-isometry, when does there exist a lifting data set $\LDS$ so that the set of all contractive
interpolants for the associated relaxed commutant lifting problem is given by
\begin{equation}\label{VtoB}
\left\{ B_V=\mat{c}{A\\\Ga_{H_V} D_A}\mid H_V=\sR_\Psi[V]\text{ for some }V\in\bS(\cE,\cE') \right\}.
\end{equation}
If $\LDS$ has this property, we say it is a lifting data set {\em associated with the Redheffer
representation $\sR_\Psi$}. Note that the lifting data set $\LDS$ should satisfy $\cD_A=\cU$ and
$\cD_{T'}=\cY$.

There is a different kind of inverse problem for interpolation and commutant lifting problems that
goes back to the work of Adamjan-Arov-Kre\u\i n \cite{AAK68,AAK71} on the Nehari problem, and has
since been considered in various settings \cite{A89,Arov90,K95a,K95b,VY,BK}. In its original form,
for the Nehari problem, this inverse problem takes the form of the question which $\gamma$-generating
pairs are also Nehari pairs. An extension to the classical commutant lifting setting is posed
and solved by J.A. Ball and A. Kheifets in \cite{BK}. For such problems the data of the considered
problem is known, and a Redheffer representation is given for which it is also known that its range
forms a subset of the set of all solutions for the problem defined by the given data. The inverse
problem is then to determine if the solutions generated by the Redheffer transformation are in fact
all solutions. In the present paper, the inverse problem is to construct a data set so that the given
Redheffer representation provides all solutions. As a consequence, the solution criteria we find here 
looks quite different from the one obtained in \cite{BK}.

Somewhat surprisingly, it turns out that for the problem considered here, it is always possible to
find an associated lifting data set.

%
%

\begin{theorem}\label{T:inverse}
Let $\sR_\Psi$ be a Redheffer representation with coefficients as a in \eqref{RedCoefs} so
that $\Psi_{1,1}(0)=0$ and  the coefficient matrix \eqref{RedMat} is a co-isometry. Then there
exists a lifting data set $\LDS$  associated with $\sR_\Psi$.  More precisely,
if $\Phi_{1,1}$, $\Phi_{1,2}$, $\Phi_{2,1}$ and $\Phi_{2,2}$ are the Redheffer
coefficients associated with $\LDS$ via \eqref{RCLcoefs}, then there exist a unitary
operator $\psi$ mapping $\cE$ onto $\cG$ and a co-isometry $\vph$ from $\cE'$ onto
$\cD_{\om^*}$ such that for each $\la\in\BD$
\begin{equation}\label{CoefsRelation}
\mat{cc}{\psi&0\\0&I_{\cY}}
\mat{cc}{\Psi_{1,1}(\la)&\Psi_{1,2}(\la)\\\Psi_{2,1}(\la)&\Psi_{2,2}(\la)}=
\mat{cc}{\Phi_{1,1}(\la)&\Phi_{1,2}(\la)\\\Phi_{2,1}(\la)&\Phi_{2,2}(\la)}\mat{cc}{\vph&0\\0&I_{\cU}}.
\end{equation}
Moreover, if in addition the coefficient matrix \eqref{RedMat} is unitary, then $\vph$ is unitary.
\end{theorem}

Since the Redheffer coefficients \eqref{RCLcoefs} are completely determined by
the underlying contraction $\om$, we see that the construction goes along the following path.
\begin{figure}[h]\label{path}
\setlength{\unitlength}{0.09in}
\centering
\begin{picture}(1,1)
\put(-26,0){$\begin{array}{c}\text{Lifting data set}\\\LDS\end{array}$}
\put(-13,0){$\boldsymbol\Longrightarrow$}
\put(-10,0){$\begin{array}{c}\text{Underlying contraction}\\\om\end{array}$}
\put(7,0){$\boldsymbol\Longrightarrow$}
\put(10,0){$\begin{array}{c}\text{Redheffer coefficients}\\
\Phi_{1,1},\ \Phi_{1,2},\ \Phi_{2,1},\ \Phi_{2,2}\end{array}$}
\end{picture}
\caption{Lifting data set to Redheffer coefficients}
\end{figure}

\noindent
It was already observed in \cite{FtHK08} that any contraction $\om$ of the form
\begin{equation}\label{genom}
\om=\mat{c}{\om_1\\\om_2}:\cF\to\mat{c}{\cY\\\cU},\quad\cF\subset\cU
\end{equation}
appears as the underlying contraction of some lifting data set. We prove Theorem
\ref{T:inverse} in Section \ref{S:inverseT} below by extracting a contraction
$\om$ of the form \eqref{genom} from the Redheffer coefficients.

Theorem \ref{T:inverse} guarantees the existence of a lifting data set with the
required properties. Moreover, as we will see in the course of the proof in Section
\ref{S:inverseT} below, the lifting data set can in fact be constructed explicitly from the
given coefficients. We already observed that for an associated lifting data set $\LDS$ we
have $\cD_A=\cU$ and $\cD_{T'}=\cY$. In addition, the space $\cF$ on which the underlying
contraction $\om$ associated with the lifting data set is defined and the unitary operator
$\psi$ can also be extracted from the given coefficients. Indeed, from the definition
of $\Phi_{1,2}$ in \eqref{RCLcoefs} and the relation \eqref{CoefsRelation} it follows that
$\cF=\kr \Psi_{1,2}(0)$, while $\psi^*=\Psi_{1,2}(0)|_{\cU\ominus\cF}$. Hence the
Hilbert spaces between which the contraction $\om$ underlying $\LDS$ acts are determined by
the Redheffer coefficients.

Note that, as one can see from \eqref{CoefsRelation}, $\Psi_{1,1}$, $\Psi_{1,2}$, $\Psi_{2,1}$
and $\Psi_{2,2}$ are, in general, not the actual coefficients associated with the constructed
data set, not even up to unitary transformations on the spaces $\cE$ and $\cE'$, unless the
operator $\vph$ is unitary. If $\vph$ is not unitary and $\cE_0'=\kr\vph$, then
$\Psi_{1,1}(\la)|_{\cE_0'}=0$ and
$\Psi_{2,1}(\la)|_{\cE_0'}=0$ for each $\la\in\BD$, which only adds to the non-uniqueness
in the Redheffer representation:
For $V,V'\in\bS(\cE,\cE')$ with $P_{\cE'\ominus\cE_0'}V(\la)=P_{\cE'\ominus\cE_0'}V'(\la)$
for all $\la\in\BD$ we have $H_V=H_{V'}$.

The two steps in the construction, lifting data set to underlying contraction and
underlying contraction to Redheffer coefficients, are in general not one-to-one.
Many lifting data sets can have the same underlying contraction and for each
admissible set of Redheffer coefficients there can be more than one underlying
contraction. However, it turns out that under some additional assumptions all
contractions $\om$ that define a fixed set of Redheffer coefficients via \eqref{RCLcoefs}
are unique up to a unitary transformation; see Proposition \ref{P:unitRedmat} below.

Let $\om'$ be the underlying contraction of another lifting data set associated with the Redheffer
coefficients \eqref{RedCoefs}, and assume that the conditions of Theorem \ref{T:inverse} are
met. Then $\om'$ also maps $\cF=\kr \Psi_{1,2}(0)$ into $\cY\oplus\cU$. We say that $\om$ and
$\om'$ are {\em unitarily equivalent} if there exists a unitary operator $\Th$ on $\cU$
for which $\cF$ is a reducing subspace ($\Th\cF\subset\cF$ and $\Th^*\cF\subset\cF$) such that
\[
\om_1\Pi_\cF\Th=\om_1'\Pi_\cF\ands \om_2\Pi_\cF\La=\La\om_2'\Pi_\cF.
\]
It is not difficult to see that unitarily equivalent underlying contractions define the same
Redheffer coefficients via \eqref{RCLcoefs}. If, in addition, the coefficient matrix \eqref{RedMat}
is unitary, then the converse is also true.

\begin{proposition}\label{P:unitRedmat}
Let $\LDS$ and $\{\wtilA,\wtilT',\wtilU',\wtilR,\wtilQ\}$ be lifting data sets associated
with Redheffer coefficients \eqref{RedCoefs} that satisfy the conditions of
Theorem~\ref{T:inverse}. Assume that the coefficient matrix \eqref{RedMat} is unitary. Then
the underlying contractions of $\LDS$ and $\{\wtilA,\wtilT',\wtilU',\wtilR,\wtilQ\}$ are
unitarily equivalent.
\end{proposition}

Since the steps in the construction are in general not one-to-one, it is of interest which
properties are invariant under these steps. Two invariants are known:
\begin{itemize}
\item[(i)] $R^*R=Q^*Q$ $\quad\Llra\quad$ $\om$ is an isometry;\\[.005cm]

\item[(ii)]$
\begin{array}{c}
\om\text{ is an isometry and}\\ \om_2\Pi_\cF\text{ is strongly stable}\end{array}$
$\ \Llra\ $
$
\begin{array}{c}
\text{the coefficient matrix}\\ \text{is unitary.}\end{array}$
\end{itemize}
If the operators in the lifting data set satisfy:
$R^*R=Q^*Q$, $\|A\|<1$ and $R$ is left invertible (or equivalently $R^*R$
is invertible), then we are in case (ii),
but these conditions on the lifting data set are only sufficient; in fact, any
$\om$ of the form \eqref{genom} is the underlying contraction of a lifting data
set with $A$ a co-isometry (see the proof of Theorem \ref{T:inverse} below).

We now return to the topic of the non-uniqueness in the Redheffer representation of
Theorem \ref{T:RCLRedheffer}. Let $B$ be a contractive interpolant for the lifting
data set $\LDS$, and define
\begin{equation}\label{nonuniqueSet}
\sV_{B}=\{V\in\bS(\cG,\cD_{\om^*}) \mid B=B_V\}.
\end{equation}
In \cite{FtHK06b} necessary and sufficient conditions were obtained that guarantee that the set
$\sV_B$ consists of one element only. To state the criteria some additional notation is required.
Assume that $B$ is given in the form \eqref{RCLsols2}, i.e., we have a contraction
$\Ga\in\sL(\cD_A,H^2(\cD_{T'}))$ such that
\[
B=\mat{c}{A\\\Ga D_A}.
\]
As in \cite{FtHK06b} we define a contraction $\om_B$ by
\begin{equation}\label{omB}
\om_B:\cF_B=\ov{D_{\Ga}\cF}\to\cD_{\Ga},\quad \om_B D_{\Ga}|_\cF=D_{\Ga}\om_2.
\end{equation}
It turns out that $\om_B$ is an isometry if and only if $\om$ is an isometry, i.e.,
if and only if $R^*R=Q^*Q$. Hence the statement ``$\om_B$ is an isometry''
holds true, or not, independent of the contractive interpolant $B$ in question.
Set $\cG_B=\cD_\Ga\ominus\cF_B$. Using this contraction $\om_B$ a one-to-one map
from the Schur class $\bS(\cG_B,\cD_{\om_B^*})$ onto the set $\sV_B$ can be
constructed; see Theorem 1.2 and the subsequent paragraph in \cite{FtHK06b}.
Since $\cG_B=\{0\}$ if and only if $\cF_B=\cD_{\Ga_H}$ and $\cD_{\om_B^*}=\{0\}$
if and only if $\om_B$ is a co-isometry, we have the following criterion for the
set $\sV_B$ to be a singleton.

\begin{proposition}\label{P:NonProper1}
Let $B$ be a contractive interpolant for $\LDS$. Let $\om_B$ be the contraction
from $\cF_{H}$ into $\cD_{\Ga_H}$ defined by \eqref{omB}. Then the set $\sV_B$
is a singleton if and only if $\cF_B=\cD_{\Ga_H}$ or $\om_B$ is a co-isometry.
\end{proposition}

Proposition \ref{P:NonProper1} gives necessary and sufficient conditions under which
the set $\sV_B$ consists of just one element for each contractive interpolant
$B$ individually. It is not clear how to derive necessary and sufficient conditions
for the map $V\mapsto B_V$ to be one-to-one from this result; sufficient conditions
obtained in the literature are: $\cF=\cD_A$, $\om$ is a co-isometry or $\om$ is an isometry
and $\ov{\om_2\cF}=\cD_A$; see Proposition 4.2.8 in \cite{tH07} and the proof of
Theorem 1.3 in \cite{FtHK06b}. However, the fact that the contraction $\om_B$ is an
isometry either for all contractive interpolants $B$ or for none gives hope
that a global result may be obtainable.

The next theorem presents a global result when we restrict to the set of contractive interpolants
$B_V$ where $V$ is from the open ball $\bS_0(\cG,\cD_{\om^*})$ of
$\bH^\infty(\cG,\cD_{\om^*})$, that is, from the set
\[
\bS_0(\cG,\cD_{\om^*}):=\{V\in\bS(\cG,\cD_{\om^*})\colon \|V\|_\infty<1\},
\]
and under the additional assumptions the underlying contraction $\om$ is an
isometry and $\om_2\Pi_\cF$ is strongly stable. The latter happens
in many of the examples considered in \cite{FFGK98,FFK02} if
the operator $A$ is a strict contraction.

\begin{theorem}\label{T:NonProper2}
Let $\LDS$ be a lifting data set such that the underlying contraction $\om$ is
an isometry and $\om_2\Pi_\cF$ on $\cD_A$ is strongly stable. Then
for each $V\in\bS_0(\cG,\cD_{\om^*})$ we have $\sV_{B_V}=\{V\}$ if and only if
\[
\Ga_{\Phi_{1,2}}\cD_A=\{\Pi_{\cG}(I-\la\om_2\Pi_\cF)u\mid u\in\cD_A\}
\]
is dense in $H^2(\cG)$. In particular, if $\sV_{B_V}=\{V\}$
for some $V\in\bS_0(\cG,\cD_{\om^*})$, then $\sV_{B_V}=\{V\}$ for all
$V\in\bS_0(\cG,\cD_{\om^*})$.
\end{theorem}

Given a lifting data set $\LDS$ that satisfies the requirement of Theorem \ref{T:NonProper2},
it follows that the Redheffer representation of Theorem \ref{T:RCLRedheffer} cannot be
one-to-one in case $\cD_{A}$ is finite dimensional, unless $\cG=\{0\}$. This occurs,
for instance, in the examples of the relaxations of classical interpolation problems
introduced in \cite{FFK02}, where the common ingredient is that $R$ and $Q$ are
operators from $\cV^{n-1}$ into $\cV^n$ of the form
\[
R=\mat{c}{I_{\cV^{n-1}}\\0}\ands
Q=\mat{c}{0\\I_{\cV^{n-1}}},
\]
assuming that $\tu{dim}\, \cV<\infty$ and $\|A\|<1$.

Another examples of a global result for the set of contractive interpolants $B_V$ associated
with $\bS_0(\cG,\cD_{\om^*})$ is the harmonic maximum principle for the classical commutant
lifting problem due to A.~Biswas \cite{B97}. As a result of our investigations of Redheffer
representations in Section \ref{S:Redheffer} below we obtain the following
extension of Theorem 3.1 from \cite{B97} to the relaxed commutant lifting setting.

\begin{theorem}\label{T:MaxPrinciple}
Given a lifting data set $\LDS$ we have:
\begin{itemize}
\item[(1)] If there exists a $\wtilV\in\bS_0(\cG,\cD_{\om^*})$ so that
$\|B_{\wtilV}\|<1$, then $\|B_V\|<1$ for all $V\in\bS_0(\cG,\cD_{\om^*})$.

\item[(2)] If there exists a $\wtilV\in\bS_0(\cG,\cD_{\om^*})$ so that
$\|B_{\wtilV}\|=1$, then $\|B_V\|=1$ for all $V\in\bS_0(\cG,\cD_{\om^*})$.

\item[(3)] Moreover, for $V\in\bS_0(\cG,\cD_{\om^*})$ and $\wtilV\in\bS(\cG,\cD_{\om^*})$
we have $\kr D_{B_V}\subset\kr D_{B_{\wtilV}}$. In particular,
$\kr D_{B_V}=\kr D_{B_{\wtilV}}$ in case $\wtilV,V\in\bS_0(\cG,\cD_{\om^*})$.
\end{itemize}
\end{theorem}

When restricting to the classical commutant lifting setting, Part 2 coincides with
Theorem 3.1 from \cite{B97}. Part 1 is just the contrapositive of Part 2. The techniques
from \cite{B97} can be used to prove Part 1 as well, but in the present paper we give a
shorter, more intuitive proof. The last part of Theorem \ref{T:MaxPrinciple} seems to be
a new result, even for classical commutant lifting. Part 3 and the fact that a contractive
interpolant $B$ with $\|B\|=\|A\|$ always exists imply the following corollary.

\begin{corollary}\label{C:nonnormattaining}
Let $\LDS$ be a lifting data set with $\|A\|<1$. Then $\kr D_{B_V}=\{0\}$
for each $V\in\bS_0(\cG,\cD_{\om^*})$, i.e., for any $h\in\cH$ and
$V\in\bS_0(\cG,\cD_{\om^*})$ we have $\|B_Vh\|<\|h\|$.
\end{corollary}

The paper consists of 6 sections, not counting the present introduction. We start with
a section where some preliminary system theory results that will be used throughout
the paper is presented. In Section \ref{S:inverseT} we prove Proposition \ref{P:unitRedmat} and
Theorem \ref{T:inverse}. The following section contains a proof of Theorem \ref{T:NonProper2}
restricted to the case that the contractive interpolant in question is the central contractive
interpolant. We then proceed with an intermezzo about general Redheffer representations in
Section \ref{S:Redheffer}.
Finally, Theorems \ref{T:NonProper2} and \ref{T:MaxPrinciple} are proved in Sections
\ref{S:nonproper} and \ref{S:maxprinciple}, respectively. The order in the proofs goes in the
reversed direction: we first prove analogous results on the level of Redheffer representations
before proving Theorems \ref{T:NonProper2} and \ref{T:MaxPrinciple}.

\setcounter{equation}{0}
\section{System theory preliminaries}\label{S:Systems}

In this section we review some results from linear system theory that will be useful
in the sequel. See \cite{F87,F81} for a general overview.

\paragraph{\bf contractive systems.}
A {\em contractive system}
is a quadruple $\{\SO,B,C,D\}$ consisting of Hilbert space operators: $\SO$ on $\cX$,
$B$ from $\cU$ to $\cX$, $C$ from $\cX$ to $\cY$ and $D$ mapping $\cU$ into $\cY$
such that the {\em system matrix}
\begin{equation}\label{SysMat}
\mat{cc}{\SO&B\\C&D}:\mat{c}{\cX\\\cU}\to\mat{c}{\cX\\\cY}
\end{equation}
is a contraction. Of particular interest will be the case when the system matrix
is a co-isometry; we then say that $\{\SO,B,C,D\}$ is a {\em co-isometric system}.
The system $\{\SO,B,C,D\}$ is said to be {\em strongly stable}
in case the state operator $\SO$ is strongly stable, that is, if for each $x\in\cX$
the sequence $\SO^nx$ converges to zero as $n\to\infty$.
Since $\SO$ is contractive, we can define analytic functions $F$ and $W$ on $\BD$ by
\begin{equation}\label{MultObs}
F(\la)=D+\la C(I-\la \SO)^{-1}B\ands
W(\la)=C(I-\la \SO)^{-1}\quad\quad(\la\in\BD).
\end{equation}
We refer to $F$ and $W$ as the {\em transfer function} and {\em observability function}
associated with $\{\SO,B,C,D\}$, respectively. {}From the fact that the system matrix
\eqref{SysMat} is a co-isometry it follows that $F\in\bS(\cU,\cY)$ and
$W\in\bH^2_\tu{ball}(\cX,\cY)$. Hence $F$ defines a contractive multiplication
operator $M_F$ from $H^2(\cU)$ to $H^2(\cY)$, and $W$ defines a contraction $\Ga_W$ from
$\cX$ into $H^2(\cY)$ via \eqref{GaH}. The system $\{\SO,B,C,D\}$ is said to be {\em observable}
in case $\kr \Ga_W=\{0\}$.

Two co-isometric systems $\{\SO_1,B_1,C_1,D_1\}$ and $\{\SO_2,B_2,C_2,D_2\}$ with the
same input spaces $\cU$ and output spaces $\cY$ but possibly different state spaces $\cX_1$
respectively $\cX_2$ are said to be {\em unitarily equivalent} if $D_1=D_2$ and
there exists a unitary operator $\Th$ mapping $\cX_1$ onto $\cX_2$ such that
\[
\Th\SO_1=\SO_2\Th,\quad C_1=C_2\Th,\quad \Th B_1=B_2.
\]
One easily verifies that unitarily equivalent systems have the same transfer function.
The converse statement is true if the co-isometric systems are also observable.

\begin{theorem}\label{T:CoOpMat}
Let $\{\SO,B,C,D\}$ be a co-isometric system with transfer function $F\in\bS(\cU,\cY)$
and observability function $W\in\bH^2_\tu{ball}(\cX,\cY)$. Then the operator
\begin{equation}\label{opmat}
\mat{cc}{\!\!M_F&\Ga_W\!\!}:\mat{c}{H^2(\cU)\\\cX}\to H^2(\cY)
\end{equation}
is a co-isometry. Moreover, the operator \eqref{opmat} is unitary if and only if the
system $\{\SO,B,C,D\}$ is strongly stable and the system matrix \eqref{SysMat} is unitary.
\end{theorem}

\begin{proof}[\bf Proof]
The statement for the case that \eqref{opmat} is unitary follows from
Theorem III.10.4 in \cite{FFGK98} and the fact that the positive operator
$\De$ on $\cX$ defined by
\[
\De^2=\tu{strong--}\!\!\lim_{n\to\infty} \SO^{*n}\SO^n
\]
is zero if and only if $\SO$ is strongly stable. For the statement on co-isometric
systems see Theorem 1.3 in \cite{tH2}.
\end{proof}

The first part of Theorem \ref{T:CoOpMat} has the following converse result.

\begin{theorem}\label{T:CoOpMatInv}
Let $F\in\bS(\cU,\cY)$ and $W\in\bH^2(\cX,\cY)$ be functions such that
the operator matrix \eqref{opmat} is a co-isometry. Then there is a co-isometric
system $\{\SO,B,C,D\}$ so that $F$ and $W$ are the associated transfer and
observability functions.
\end{theorem}

\begin{proof}[\bf Proof]
It follows from \cite{Ando90} that the function $F$ is the transfer
function of an observable co-isometric system $\{\wtilSO,\wtilB,\wtilC,\wtilD\}$, say
with state space $\widetilde\cX$. {}From Theorem
\ref{T:CoOpMat} we know that the operator matrix
$\mat{cc}{M_{F}&\Ga_{\wtilW}}$ is a co-isometry. Hence
\[
M_FM_{F^*}+\Ga_{\wtilW}\Ga_{\wtilW}^*=I=M_FM_{F^*}+\Ga_{W}\Ga_{W}^*,\quad\text{ and thus }\quad
\Ga_{\wtilW}\Ga_{\wtilW}^*=\Ga_{W}\Ga_{W}^*.
\]
By Douglas factorization lemma \cite{D66} there exists a
unitary operator $\La$ mapping $\ov{\im\Ga_{\wtilW}^*}=\widetilde\cX$ onto
$\cN=\ov{\im\Ga_{W}^*}$ defined by the identity $\Ga_{\wtilW}\La^*=\Ga_{W}$.
Now set $\cM=\cX\ominus\cN$ and define $\{\SO,B,C,D\}$ to be the co-isometric
system given by the system matrix
\[
\mat{cc}{\SO&B\\C&D}=\mat{cc|c}{\La&0&0\\0&I_{\cM}&0\\\hline0&0&I_\cY}
\mat{cc|c}{\wtilSO&0&\wtilB\\0&I_\cM&0\\\hline\wtilC&0&\wtilD}
\mat{cc|c}{\La^*&0&0\\0&I_{\cM}&0\\\hline0&0&I_\cY}.
\]
Then $\{\SO,B,C,D\}$ and $\{\wtilSO,\wtilB,\wtilC,\wtilD\}$ are unitarily equivalent,
and thus $F$ is also the transfer function of $\{\SO,B,C,D\}$. Moreover, the
observability function $\widehat W$ of $\{\SO,B,C,D\}$ satisfies $\widehat
W(\la)=\mat{cc}{\wtilW(\la)\La^*&0}$ for each $\la\in\BD$. The
latter implies that
\[
\Ga_{\widehat W}=\mat{cc}{\Ga_{\wtilW}\La^*&0}=\Ga_W.
\]
But then $\widehat W=W$, which proves our claim.
\end{proof}

\paragraph{\bf Redheffer cascading systems.}
Suppose we are given Hilbert spaces $\cU_1$, $\cU_2$, $\cY_1$, $\cY_2$, $\cX$
and $\cX'$ and two operators $M_1$ mapping $\cX\oplus\cU_1$ into $\cX'\oplus\cY_1$
and $M_2$ mapping $\cX'\oplus\cU_2$ into $\cX\oplus\cY_2$ with operator matrix
decompositions
\begin{equation}\label{M1M2}
M_1=\mat{cc}{\SO_1&B_1\\C_1&D_1}:\mat{c}{\cX\\\cU_1}\to\mat{c}{\cX'\\\cY_1},\ \
M_2=\mat{cc}{\SO_2&B_2\\C_2&D_2}:\mat{c}{\cX'\\\cU_2}\to\mat{c}{\cX\\\cY_2}.
\end{equation}
Furthermore, assume that $I_\cX-\SO_2\SO_1$ is invertible. Then
$I_{\cX'}-\SO_1\SO_2$ is invertible, and we can form the {\em Redheffer product}
$M_1\circ M_2$ of $M_1$ and $M_2$, which is the operator from $\cU_1\oplus\cU_2$
into $\cY_1\oplus\cY_2$ given by the operator matrix decomposition
\begin{equation}\label{FeedbackMat}
\mat{cc}{D_1+C_1\SO_2(I_{\cX'}-\SO_1\SO_2)^{-1}B_1&C_1(I_{\cX}-\SO_2\SO_1)^{-1}B_2\\
C_1(I_{\cX'}-\SO_1\SO_2)^{-1}B_1&D_2+C_2\SO_1(I_{\cX}-\SO_2\SO_1)^{-1}B_2}.
\end{equation}
For an elaborate discussion on Redheffer products we refer to \cite[Chapter XIV]{FF90}.
The next proposition shows how norm properties of $M_1$ and $M_2$ carry over to
$M_1\circ M_2$.

\begin{proposition}\label{P:CoOpMatInv}
In case the block operator matrices $M_1$ and $M_2$ in \eqref{M1M2} are both
contractive, isometric, co-isometric or unitary, and $I-X_2X_1$ is invertible,
the Redheffer product $M_1\circ M_2$ is also contractive, isometric,
co-isometric or unitary, respectively.
\end{proposition}

\begin{proof}[\bf Proof]
Since we assume $I_\cX-\SO_2\SO_1$ to be invertible, the case where $M_1$ and $M_2$
are both contractive follows from Lemma XIV.1.2 in \cite{FF90}. The case where
$M_1$ and $M_2$ are isometries follows from Lemma XIV.1.2 in \cite{FF90} as well, using
identity (1.15) on Page 433 of \cite{FF90}. The co-isometric case is obtained by
applying the statement for the isometric case to $M_1^*$ and $M_2^*$, and finally
the unitary case follows from the result when $M_1$ and $M_2$ are both isometric
and co-isometric.
\end{proof}

\setcounter{equation}{0}
\section{Proofs of Proposition  \ref{P:unitRedmat}
and Theorem \ref{T:inverse}}\label{S:inverseT}

We start this section with a proof of Theorem \ref{T:inverse},
followed by a proof of Proposition \ref{P:unitRedmat}.

\begin{proof}[\bf Proof of Theorem \ref{T:inverse}]
Assume that $\cU$, $\cY$, $\cE$ and $\cE'$ are Hilbert spaces, and
\[
\Psi_{1,1}\in\bS(\cE',\cE),\ \
\Psi_{1,2}\in\bH^2_\tu{ball}(\cU,\cE),\ \
\Psi_{2,1}\in\bS(\cE',\cY),\ \
\Psi_{2,2}\in\bH^2_\tu{ball}(\cU,\cY),
\]
with $\Psi_{1,1}(0)=0$ such that the coefficient matrix
\eqref{RedMat} is a co-isometry. Set
\[
F(\la)=\mat{c}{\Psi_{1,1}(\la)\\\Psi_{2,1}(\la)}\in\bS(\cE',\cE\oplus\cY)\ands
W(\la)=\mat{c}{\Psi_{1,2}(\la)\\\Psi_{2,2}(\la)}\in\bH^2_\tu{ball}(\cU,\cE\oplus\cY).
\]
Then $\mat{cc}{M_F&\Ga_W}$ is a co-isometry from $H^2(\cE')\oplus\cU$ into $H^2(\cE\oplus\cY)$.
According to Theorem \ref{T:CoOpMatInv} there exists a co-isometric system
$\{\SO,B,\sbm{C_1\\C_2},\sbm{D_1\\D_2}\}$ with system matrix
\[
N=\mat{cc}{\SO&B\\C_1&D_1\\C_2&D_2}:\mat{c}{\cU\\\cE'}\to\mat{c}{\cU\\\cE\\\cY}
\]
such that $F$ and $W$ are the corresponding transfer function and observability function, respectively.
The fact that $\Psi_{1,1}(0)=0$ implies that $D_1=0$, and thus, since the system matrix $N$ is a
co-isometry, that $C_1=\Psi_{1,2}(0)$ is a co-isometry. Set $\cF=\kr C_1$ and $\cG=\cU\ominus\cF$. Then $C_1$
maps $\cG$ isometrically onto $\cE'$. Define $\om_1=C_2|_\cF$ and $\om_2=\SO|_\cF$, i.e.,
\[
N=\mat{cc}{\om_2\Pi_\cF&B\\C_1P_\cG&0\\\om_1\Pi_\cF&D_1}.
\]
By a special case of Parrot's lemma \cite{P78}, again using that
$N$ is a co-isometry, there exists a
co-isometry $\vph$ mapping $\cE'$ onto $\cD_{\om^*}$
such that
\begin{equation}\label{Defectom*}
\mat{cc}{B\\D_1}=D_{\om^*}\vph,\quad\text{where}\quad
\om=\mat{c}{\om_1\\\om_2}:\cF\to\mat{cc}{\cY\\\cU}.
\end{equation}
Let $\psi$ be the unitary map from $\cG$ onto $\cE$ given by $\psi=C_1|_\cG$. Then
\[
\mat{ccc}{I_\cU&0&0\\0&\psi&0\\0&0&I_\cY}N\mat{cc}{I_\cU&0\\0&\vph^*}
=\mat{cc}{\om_2\Pi_\cF&\Pi_\cU D_{\om^*}\\\Pi_\cG&0\\\om_1\Pi_\cF&\Pi_\cY D_{\om^*}}.
\]
Thus $\{\SO,B,\sbm{C_1\\C_2},\sbm{D_1\\D_2}\}$ is unitarily equivalent to the system
$\{\widehat{\SO},\widehat{B},\widehat{C},\widehat{D}\}$ given by
\begin{equation}\label{REDsysmat}
\mat{cc}{\widehat{\SO}&\widehat{B}\\\widehat{C}&\widehat{D}}
=\mat{c|c}{\om_2\Pi_\cF&\Pi_\cU D_{\om^*}\\\hline\Pi_\cG&0\\\om_1\Pi_\cF&\Pi_\cY D_{\om^*}}
:\mat{c}{\cU\\\cD_{\om^*}}\to\mat{c}{\cU\\\cG\\\cY}.
\end{equation}
One easily verifies that
\[
\mat{c}{\Phi_{1,1}(\la)\\\Phi_{2,1}(\la)}
=\widehat{D}+\la\widehat{C}(I-\la \widehat{\SO})^{-1}\widehat{B}\ands
\mat{c}{\Phi_{1,2}(\la)\\\Phi_{2,2}(\la)}=\widehat{C}(I-\la \widehat{\SO})^{-1},
\]
where $\Phi_{1,1}$ and $\Phi_{2,1}$ are defined by \eqref{RCLcoefs} with $\om$
as in \eqref{Defectom*}. A straightforward computation shows that the
transfer function $\widehat{F}$ and observability function $\widehat{W}$ of the co-isometric system
$\{\widehat\SO,\widehat{B},\widehat{C},\widehat{D}\}$ given in \eqref{REDsysmat}
relate to $F$ and $W$ via
\[
\mat{cc}{\psi&0\\0&I_\cY}\mat{cc}{F(\la)&W(\la)}
=\mat{cc}{\widehat{F}(\la)&\widehat{W}(\la)}\mat{cc}{\vph&0\\0&I_\cU}\quad(\la\in\BD).
\]
Thus, to prove Theorem \ref{T:inverse} it remain to show that the contraction $\om$
appear as the underlying contraction of a lifting data set. However, it is known from
\cite{FtHK08} that any contraction $\om$ of the form in \eqref{Defectom*} with $\cF\subset\cU$, is
the underlying contraction of some lifting data set. In fact, one easily verifies that
the operators
\[
A=\mat{cc}{I_\cY&0}:\mat{c}{\cY\\\cU}\to\cY,\ T'=0\ons\cY,\ R=\om,Q=\Pi_\cF^*:\cF\to\cY\oplus\cU
\]
define a lifting data set $\LDS$, with $U'$ the Sz.-Nagy-Sch\"affer isometric lifting
of $T'$, which has $\om$ as its underlying contraction.

To complete the proof, assume that the coefficient matrix $K_0$ in \eqref{RedMat} is unitary.
Then $K_0|_{H^2(\cE')}$ is an isometry, and thus, in particular, for each $e\in\cE'$
there exists a $\la\in\BD$ such that $\Psi_{1,1}(\la)e\not=0$ or $\Psi_{2,1}(\la)e\not=0$.
So, by \eqref{CoefsRelation} and the fact that $\psi$ is unitary,
$\vph e\not=0$ for all $e\in\cE'$, and thus $\vph$ is, in fact, unitary.
\end{proof}

\begin{proof}[\bf Proof of Proposition \ref{P:unitRedmat}]
Let $\om$ and $\om'$ denote the underlying contractions of the data sets $\LDS$ and
$\{\wtilA,\wtilT',\wtilU',\wtilR,\wtilQ\}$, respectively.

Since the coefficient matrix \eqref{RedMat} is assumed to be unitary, we may without loss
of generality assume that the Redheffer coefficients defined by $\om$ and $\om'$ via
\eqref{RCLcoefs} are equal to the given Redheffer coefficients \eqref{RedCoefs}. Indeed,
this is the case because, according to the last part of Theorem \ref{T:inverse}, the three
sets of Redheffer coefficients are the same up to unitary transformations between
$\cE$ and $\cG=\cU\ominus\cF$ and between $\cE'$ and $\cD_{\om^*}$ and $\cD_{\om^{'*}}$,
and these unitary transformations do not affect whether $\om$ and $\om'$ are unitarily
equivalent, or not.

Now let $\{\widehat{\SO},\widehat{B},\widehat{C},\widehat{D}\}$ be the co-isometric
system given by \eqref{REDsysmat} and define
$\{\widehat{\SO}',\widehat{B}',\widehat{C}',\widehat{D}'\}$ accordingly with $\om$
replaced by $\om'$. Then both systems have the same transfer function $F$ and observability
function $W$ that are given by
\[
F(\la)=\mat{c}{\Psi_{1,1}(\la)\\\Psi_{2,1}(\la)}
\ands
W(\la)=\mat{c}{\Psi_{1,2}(\la)\\\Psi_{2,2}(\la)}.
\]
By assumption, the coefficient matrix is unitary. This implies that the operator $\Ga_W$,
defined via \eqref{GaH}, is an isometry. In particular, $\kr \Ga_W=\{0\}$. Hence the two
co-isometric systems $\{\widehat{\SO},\widehat{B},\widehat{C},\widehat{D}\}$ and
$\{\widehat{\SO}',\widehat{B}',\widehat{C}',\widehat{D}'\}$ are observable and have the
same transfer function. This implies that they are unitarily equivalent. It then follows
from the definition of the systems $\{\widehat{\SO},\widehat{B},\widehat{C},\widehat{D}\}$ and
$\{\widehat{\SO}',\widehat{B}',\widehat{C}',\widehat{D}'\}$ that the operator $\Th$ on $\cU$
that establishes this unitary equivalence is also the operator that establishes the unitary
equivalence of $\om$ and $\om'$.
\end{proof}

\setcounter{equation}{0}
\section{Proof of Theorem \ref{T:NonProper2} for the central solution}

Let $B_c$ denote the central contractive interpolant, that is, $B_c=B_{\CS}$, where
$\CS$ denotes the constant Schur class function whose value is the zero operator
in $\sL(\cG,\cD_{\om^*})$. In this section we prove Theorem \ref{T:NonProper2} restricted
to the case that the contractive interpolant $B$ is equal to $B_c$. To be more precise,
we prove the following proposition.

\begin{proposition}\label{P:CentralNonProper}
Let $\LDS$ be a lifting data set such that the underlying contraction $\om$ is
an isometry and $\om_2\Pi_\cF$ on $\cD_A$ is strongly stable.
Then the set $\sV_{B_c}$ in \eqref{nonuniqueSet} consists of just one element
if and only if $\Ga_{\Phi_{1,2}}\cD_A$ is dense in $H^2(\cG)$.
\end{proposition}

In order to prove Proposition \ref{P:CentralNonProper} it is convenient to first prove the next lemma.
Recall that $E_\cG\in\sL(\cG, H^2(\cG))$ stands for the canonical embedding of $\cG$ into the set
of constant functions in $H^2(\cG)$.

\begin{lemma}\label{L:DenseRange}
Let $\cK\subset H^2(\cG)$ be an invariant subspace of the backward shift
$S_\cG^*$. Then $\cG_0:=E_\cG\cG\subset\cK$ and $S_{\cG}^*\cK=\cK$ if and only
if $\cK=H^2(\cG)$.
\end{lemma}

\begin{proof}[\bf Proof]
If $\cK=H^2(\cG)$, then obviously $\cG_0\subset\cK$ and
$S_{\cG}^*\cK=\cK$. Conversely, assume that  $\cG_0\subset\cK$ and
$S_{\cG}^*\cK=\cK$. Then $S_{\cG}^n\cG_0\subset
S_{\cG}^n\cK=S_{\cG}^nS_{\cG}^{*n}\cK$ for $n=0,1,2,\ldots$. Since
$\cG_0=S_{\cG}^0\cG_0\subset\cK$, we obtain that
$S_\cG\cG_0\subset S_\cG S_\cG^*\cK=\cK\ominus\cG_0\subset\cK$, and
with a similar argument it follows, recursively, that
\[
S_{\cG}^n\cG_0\subset S_\cG^nS_\cG^{*n}\cK=\cK\ominus(\oplus_{k=0}^{n-1}S_{\cG}^k\cG_0)\subset\cK
\]
for each nonnegative integer $n$. The latter inclusions imply that $\cK=H^2(\cG)$.
\end{proof}

\begin{proof}[\bf Proof of Proposition \ref{P:CentralNonProper}]
Note that $B_c$ is of the form \eqref{RCLsols2} with $H_V$ equal to $H_c:=\Phi_{2,2}$.
Let $\om_{B_c}$ be the contraction defined according to \eqref{omB}. It follows from
the last statement of Theorem \ref{T:RCLRedheffer} and the remark in the sentence
after \eqref{omB} that $\om_{B_c}$ is an isometry. So, by Proposition \ref{P:NonProper1},
it suffices to prove that $\Ga_{\Phi_{1,2}}\cD_A$ is dense in $H^2(\cG)$ if and only if
$\cF_{B_c}=\cD_{\Ga_{\Phi_{2,2}}}$ or $\om_{B_c}\cF_{B_c}=\cD_{\Ga_{\Phi_{2,2}}}$.

Since, by Theorem \ref{T:RCLRedheffer}, the coefficient matrix $K_0$ is unitary we have
\[
D_{\Ga_{H_c}}^2=D_{\Ga_{\Phi_{2,2}}}^2=I-\Ga_{\Phi_{2,2}}^*\Ga_{\Phi_{2,2}}
=\Ga_{\Phi_{1,2}}^*\Ga_{\Phi_{1,2}}.
\]
Hence there exists a unique unitary operator $\rho$ mapping $\cD_{H_c}$ onto
$\cK:=\ov{\Ga_{\Phi_{1,2}}\cD_A}$ defined by $\rho D_{\Ga_{\Phi_{2,2}}}=\Ga_{\Phi_{1,2}}$.
The fact that $\Ga_{\Phi_{1,2}}$ is contractive and $\Phi_{1,2}(0)=\Pi_\cG$ implies that
$\cG_0:=E_\cG\cG\subset\cK$ and that $\rho$ maps $\cF_{B_c}=D_{\Ga_{\Phi_{2,2}}}\cF$
onto $\cFt_{B_c}:=\cK\ominus\cG_0$. In particular, $\cF_{B_c}=\cD_{\Ga_{\Psi_{2,2}}}$
holds if and only if $\cG=\{0\}$, from which the denseness of $\Ga_{\Phi_{1,2}}\cD_A$
in $H^2(\cG)$ obviously follows.

Let $\tilom_{B_c}$ be the isometry from $\cFt_{B_c}$ to $\cK$
defined by $\tilom_{B_c}\!=\!\rho\om_{B_c}\rho^*|_{\cFt_{B_c}}$. Note
that
$\tilom_{B_c}\Pi_{\cFt_{B_c}}=\rho\om_{B_c}\Pi_{\cF_{B_c}}\rho^*$.
We claim that $\cK$ is invariant under $S_{\cG}^*$ and
$\tilom_{B_c}\Pi_{\cFt_{B_c}}=S_\cG^*|_\cK$. The latter, which
implies the former, follows since
\begin{eqnarray*}
\tilom_{B_c}\Pi_{\cFt_{B_c}}\Ga_{\Phi_{1,2}}
&=&\rho\om_{B_c}\Pi_{\cF_{B_c}}\rho^*\Ga_{\Phi_{1,2}}
=\rho\om_{B_c}\Pi_{\cF_{B_c}}D_{\Ga_{\Phi_{2,2}}}
=\rho\om_{B_c}D_{\Ga_{\Phi_{2,2}}}P_{\cF}\\
&=&\rho D_{\Ga_{\Phi_{2,2}}}\om_2\Pi_{\cF}
=\Ga_{\Phi_{1,2}}\om_2\Pi_{\cF}
=S^*_\cG\Ga_{\Phi_{1,2}}P_{\cF}
=S^*_\cG\Ga_{\Phi_{1,2}}.
\end{eqnarray*}
The third identity follows from the fact that $\Ga_{\Phi_{2,2}}|_\cG=0$, and thus
$D_{\Ga_{\Phi_{2,2}}}\cG=\cG$. For the last identity it is again used that $\Ga_{\Phi_{1,2}}$
is a contraction and $\Phi_{1,1}(0)=\Pi_\cG$, which implies that
$\Ga_{\Phi_{1,2}}\cG=\cG_0=\kr S^*_{\cG}$.

We have thus shown that $\cK$ is a backward shift invariant subspace of $H^2(\cG)$ with
$\cG_0\subset\cK$. It then follows from Lemma \ref{L:DenseRange} that
$\cK=H^2(\cG)$ if and only if $\tilom_{B_c}\cFt_{B_c}=S^*_{\cG}\cK=\cK$.
This is the same as saying that $\Ga_{\Phi_{1,2}}\cD_A$  is dense in
$H^2(\cG)$ if and only if $\om_{B_c}\cF_{B_c}=\cD_{\Ga_{\Phi_{2,2}}}$.
\end{proof}

As a result of Theorem \ref{T:inverse} we obtain the following analogue of
Proposition \ref{P:CentralNonProper} in the setting of general Redheffer
representations.

\begin{proposition}\label{P:CentralNonProper2}
Let $\sR_\Psi$ be a Redheffer representation with Redheffer coefficients
$\Psi_{1,1}$, $\Psi_{1,2}$, $\Psi_{2,1}$ and $\Psi_{2,2}$ as in \eqref{RedCoefs}
such that the coefficient matrix \eqref{RedMat} is unitary. Then the set
\[
\sW_{H_{\CS}}:=\{V\in\bS(\cE,\cE')\mid H_{\CS}=H_V\}
\]
is equal to $\{\CS\}$ if and only if $\Ga_{\Psi_{1,2}}\cU$ is dense in $H^2(\cE)$.
\end{proposition}

\begin{proof}[\bf Proof]
The result follows immediately from Theorem \ref{T:inverse} and
Proposition \ref{P:CentralNonProper}.
\end{proof}

\setcounter{equation}{0}
\section{Intermezzo:~Redheffer representations}\label{S:Redheffer}

Assume we are given Hilbert spaces $\cU$, $\cY$, $\cE$ and $\cE'$ and operator-valued functions
\begin{equation}\label{RedCoefs2}
\Psi_{1,1}\in\bS(\cE',\cE),\ \
\Psi_{1,2}\in\bH^2_\tu{ball}(\cU,\cE),\ \
\Psi_{2,1}\in\bS(\cE',\cY),\ \
\Psi_{2,2}\in\bH^2_\tu{ball}(\cU,\cY)
\end{equation}
with $\Psi_{1,1}(0)=0$. In the present section we derive some results for the Redheffer representation
$\sR_\Psi$ with coefficients $\Psi_{1,1}$, $\Psi_{1,2}$, $\Psi_{2,1}$ and $\Psi_{2,2}$.

In case the coefficient matrix
\begin{equation}\label{RedMat2}
K_0=\mat{cc}{M_{\Psi_{1,1}}&\Ga_{\Psi_{1,2}}\\M_{\Psi_{2,1}}&\Ga_{\Psi_{2,2}}}
\mat{c}{H^2(\cE')\\\cU}\to\mat{c}{H^2(\cE)\\H^2(\cY)}
\end{equation}
is a co-isometry, it follows from Theorems \ref{T:inverse} and \ref{T:RCLRedheffer} that
$\sR_\Psi$ maps $\bS(\cE,\cE')$ into $\bH^2_\tu{ball}(\cU,\cY)$. The next proposition
shows that this is still true if the coefficient matrix is a contraction.

\begin{proposition}\label{P:Redheffer}
Let $\Psi_{1,1}$, $\Psi_{1,2}$, $\Psi_{2,1}$ and $\Psi_{2,2}$ be as in
\eqref{RedCoefs2} with $\Psi_{1,1}(0)=0$ and assume that \eqref{RedMat2} is contractive.
Then $\sR_\Psi$ maps $\bS(\cE,\cE')$ into $\bH^2_\tu{ball}(\cU,\cY)$.
\end{proposition}

\begin{proof}[\bf Proof]
For the case that the coefficient matrix $K_0$ is a co-isometry, the
statement follows right away from Theorems \ref{T:RCLRedheffer} and \ref{T:inverse}.
Now assume that $K_0$ is contractive. By Parrot's lemma \cite{P78} there exist a
Hilbert space $\cV$ and functions $\Psi_{1,3}\in\bH^2_\tu{ball}(\cV,\cE)$ and
$\Psi_{2,3}\in\bH^2_\tu{ball}(\cV,\cY)$ such that
\[
\wtilK_0=\mat{ccc}{M_{\Psi_{1,1}}&\Ga_{\Psi_{1,2}}&\Ga_{\Psi_{1,3}}\\
M_{\Psi_{2,1}}&\Ga_{\Psi_{2,2}}&\Ga_{\Psi_{2,3}}}:
\mat{c}{H^2(\cE')\\\cU\\\cV}\to\mat{c}{H^2(\cE)\\H^2(\cY)}
\]
is a co-isometry. Set
$\widetilde\Psi_{1,2}(\la)=\mat{cc}{\Psi_{1,2}(\la)&\Psi_{1,3}(\la)}\in\bH^2_\tu{ball}(\cU\oplus\cV,\cE)$
and
$\widetilde\Psi_{2,2}(\la)=\mat{cc}{\Psi_{2,2}(\la)&\Psi_{2,3}(\la)}\in\bH^2_\tu{ball}(\cU\oplus\cV,\cY)$.
It follows that the Redheffer representation $\sR_{\widetilde\Psi}$ with Redheffer coefficients
$\Psi_{1,1}$, $\widetilde\Psi_{1,2}$, $\Psi_{2,1}$ and $\widetilde\Psi_{2,2}$ has a co-isometric
coefficient matrix, and thus maps $\bS(\cE,\cE')$ into $\bH^2_\tu{ball}(\cU\oplus\cV,\cY)$.
The statement now follows because $\sR_\Psi[V](\la)=\sR_{\widetilde\Psi}[V](\la)|_{\cU}$
for each $\la\in\BD$ and $V\in\bS(\cE,\cE')$. Hence $H_V\in\bH^2_\tu{ball}(\cU,\cY)$.
\end{proof}

As in the relaxed commutant lifting setting, the map $\sR_\Psi$ is
in general not one-to-one. This is the underlying obstruction preventing the obviously
sufficient condition ``$\cE=\{0\}$ or $\cE'=\{0\}$'' for $\sR_\Psi[\bS(\cE,\cE')]$
to be a singleton from also being necessary. Necessary and sufficient conditions under which
the range of $\sR_\Psi$ is a singleton, in the relaxed commutant lifting setting,
were obtained in \cite{LT06} for the case that $R^*R=Q^*Q$ and in \cite{tH2} for the
general case. An adaptation of the proof in \cite{tH2} gives the following result.

\begin{proposition}\label{P:UniqueSolution}
Let $\Psi_{1,1}$, $\Psi_{1,2}$, $\Psi_{2,1}$ and $\Psi_{2,2}$ be as in
\eqref{RedCoefs2} with $\Psi_{1,1}(0)=0$ and assume that the coefficient
matrix \eqref{RedMat2} is contractive. Then the range of $\sR_\Psi$
is a singleton if and only if $\Psi_{1,2}(\la)=0$ for each $\la\in\BD$ or
$\Psi_{2,1}(\la)=0$ for each $\la\in\BD$. In case the coefficient matrix
is a co-isometry, the range of $\sR_\Psi$ is a singleton if and only
if $\cE=\{0\}$ or $\Ga_{\Psi_{2,2}}$ is a co-isometry.
\end{proposition}

\begin{proof}[\bf Proof]
The statement for the case that the coefficient matrix is a co-isometry
follows from Theorem \ref{T:inverse} above and Theorem 0.2 in \cite{tH2}. The
result for the case that the coefficient matrix is contractive can
easily be deduced from the proof of Theorem~0.2 in \cite{tH2}.
\end{proof}

We now restrict our attention to the Schur class functions in the
open unit ball $\bS_0(\cE,\cE')$ of $\bH^\infty(\cE,\cE')$, i.e.,
$V\in\bH^\infty(\cE,\cE')$ with $\|V\|_\infty<1$. For
$V\in\bS_0(\cE,\cE')$ the multiplication operator $M_V$ is a
strict contraction, which enables us to express the operator $\Ga_{H_V}$
associated with $H_V=\sR_\Psi[V]$ via \eqref{GaH} explicitly as
\begin{equation}\label{GaHVexplicit}
\Ga_{H_V}=\Ga_{\Psi_{2,2}}+M_{\Psi_{2,1}}M_V(I-M_{\Psi_{1,1}}M_V)^{-1}\Ga_{\Psi_{1,2}}.
\end{equation}
We write $R_V$ for the rotation operator associated with the strict contraction $M_V$:
\begin{equation}\label{RotOperV2}
R_V:=\mat{cc}{M_V&D_{M_V^*}\\-D_{M_V}&M_V^*}
:\mat{c}{H^2(\cE)\\H^2(\cE')}\to\mat{c}{H^2(\cE')\\H^2(\cE)}.
\end{equation}
Then $R_V$ is unitary, and, since $I-M_{\Psi_{1,1}}M_V$ is invertible, we can
form the Redheffer product $K_V=R_V\circ K_0$ of $R_V$ with the coefficient
matrix $K_0$ in \eqref{RedMat2},  which gives us the following block operator
matrix
\begin{equation}\label{KV2}
K_V=\mat{cc}{M_V^*+D_{M_V}(I-M_{\Psi_{1,1}}M_V)^{-1}M_{\Psi_{1,1}}&
-D_{M_V}(I-M_{\Psi_{1,1}}M_V)^{-1}\Ga_{\Psi_{1,2}}\\
M_{\Psi_{2,1}}(I-M_VM_{\Psi_{1,1}})^{-1}D_{M_V^*}&
\Ga_{H_V}}
\end{equation}
that maps $H^2(\cE')\oplus\cU$ into $H^2(\cE)\oplus H^2(\cY)$. Note that in case
$V=\underline{0}$, the constant function with value $0\in\sL(\cE,\cE')$, $K_V$ is
equal to the coefficient matrix $K_0$.

\begin{proposition}\label{P:KVMat2}
For each $V\in\bS_0(\cE,\cE')$ the operator matrix $K_V$ is contractive,
co-isometric, isometric or unitary whenever the coefficient matrix $K_0$
in \eqref{RedMat2} is contractive, co-isometric, isometric or unitary,
respectively.
\end{proposition}

\begin{proof}[\bf Proof]
The result follows right away from Proposition \ref{P:CoOpMatInv} and
the fact that the rotation operator $R_V$ is unitary.
\end{proof}

\setcounter{equation}{0}
\section{Proof of Theorem \ref{T:NonProper2}}\label{S:nonproper}

Before giving the proof of Theorem \ref{T:NonProper2}, we first prove the
following Redheffer representation analogue.

\begin{theorem}\label{T:NonProper2RED}
Let $\Psi_{1,1}$, $\Psi_{1,2}$, $\Psi_{2,1}$ and $\Psi_{2,2}$ be as in
\eqref{RedCoefs2} with $\Psi_{1,1}(0)=0$ and assume that the coefficient
matrix $K_0$ in \eqref{RedMat2} is unitary. For each $V\in\bS_0(\cE,\cE')$ the set
\begin{equation}\label{nonuniqueSet2}
\sW_{H_V}=\{\tilV\in\bS(\cE,\cE') \mid H_V=H_{\tilV}\}
\end{equation}
is equal to $\{V\}$ if and only if $\Ga_{\Psi_{1,2}}\cU$ is dense in $H^2(\cE)$.
In particular, if $\sW_{H_V}=\{V\}$ for some $V\in\bS_0(\cE,\cE')$, then
$\sW_{H_V}=\{V\}$ for all $V\in\bS_0(\cE,\cE')$.
\end{theorem}

\begin{proof}[\bf Proof]
By Theorem \ref{T:inverse} there exist a lifting data set $\LDS$ with underlying
contraction $\om$ as in \eqref{om} and associated Redheffer coefficients
$\Phi_{1,1}$, $\Phi_{1,2}$, $\Phi_{2,1}$ and $\Phi_{2,2}$ given by \eqref{RCLcoefs},
and unitary operators $\psi$ mapping $\cE$ onto $\cG=\cD_A\ominus\cF$ and
$\vph$ from $\cE'$ onto $\cD_{\om^*}$ such that \eqref{CoefsRelation} holds.
Thus $\Psi_{1,1}$, $\Psi_{1,2}$, $\Psi_{2,1}$ and $\Psi_{2,2}$ can be identified
with coefficients $\Phi_{1,1}$, $\Phi_{1,2}$, $\Phi_{2,1}$ and $\Phi_{2,2}$, up to the
two unitary maps $\psi$ and $\vph$. In particular, the coefficient matrix \eqref{CoefMat}
is unitary as well, which, according to Theorem \ref{T:RCLRedheffer}, implies that
$\om$ is an isometry.

Now fix a $V\in\bS_0(\cE,\cE')$. The relation between the two sets of Redheffer coefficients
and the lifting data set $\LDS$ implies that we can define a contraction $\om_{H_V}$ by
\[
\om_{H_V}:\cF_{H_V}=\ov{D_{\Ga_{H_V}}\cF}\to\cD_{\Ga_{H_V}},\quad
\om_{H_V}D_{\Ga_{H_V}}|_\cF=D_{\Ga_{H_V}}\om_2.
\]
In fact, $\om_{H_V}$ is an isometry, because $\om$ is an isometry; see the paragraph
preceding Proposition \ref{P:NonProper1}. Moreover, by Proposition \ref{P:NonProper1} and from the
fact that $\om_{H_V}$ is an isometry, it follows that $\sW_{H_V}=\{V\}$ if and only if
$\cF_{H_V}=\cD_{\Ga_{H_V}}$ or $\om_{H_V}\cF_{H_V}=\cD_{\Ga_{H_V}}$.

Since the coefficient matrix \eqref{RedMat2} is unitary, we obtain from Proposition
\ref{P:KVMat2} that the operator $K_V$ in \eqref{KV2} is unitary. In particular,
\[
D_{\Ga_{H_V}}^2=I-\Ga_{H_V}^*\Ga_{H_V}=C_V^*C_V,
\]
where we define $C_V$ by
\[
C_V=Y_V\Ga_{\Psi_{1,2}}\quad\text{with}\quad Y_V=D_{M_V}(I-M_{\Psi_{1,1}}M_V)^{-1}.
\]
Set $\cK_V=\ov{C_V\cU}$. According to Douglas' factorization lemma \cite{D66},
the identity $\rho D_{\Ga_{H_V}}=C_V$ defines a unitary map $\rho$ from $\cD_{\Ga_{H_V}}$
onto $\cK_V$. Set
\[
\tilom_{H_V}=\rho\om_{H_V}\rho^*|_{\cFt_{H_V}},\quad\text{where}\quad\cFt_{H_V}=\rho\cF_{H_V}.
\]
Then $\tilom_{H_V}$ is an isometry from $\cFt_{H_V}$
into $\cK_V$, and $\sW_{H_V}=\{V\}$ if and only if
\begin{equation}\label{aaa}
\cFt_{H_V}=\ov{Y_V\Ga_{\Psi_{1,2}}\cF}=\cK_V\quad\text{or}\quad
\tilom_{H_V}\cFt_{H_V}=\ov{Y_V\Ga_{\Psi_{1,2}}\om_2\cF}=\cK_V.
\end{equation}
Since $\|M_V\|=\|V\|_\infty<1$, $D_{M_V}$ is invertible, and hence $Y_V$ is invertible.
But then \eqref{aaa} holds if and only if
\begin{equation}\label{bbb}
\ov{\Ga_{\Psi_{1,2}}\cF}=\ov{\Ga_{\Psi_{1,2}}\cU}\quad\text{or}\quad
\ov{\Ga_{\Psi_{1,2}}\om_2\cF}=\ov{\Ga_{\Psi_{1,2}}\cU}.
\end{equation}
The conditions in \eqref{bbb} are equal to those in \eqref{aaa} in case  $V=\CS$,
the constant function in $\bS(\cE,\cE')$ whose value is the zero operator. This case
was resolved in Proposition \ref{P:CentralNonProper2}, from which we infer that
$\sW_{H_V}=\{V\}$ if and only if $\Ga_{\Psi_{1,2}}\cU$ is dense in $H^2(\cE)$.
\end{proof}

\begin{proof}[\bf Proof of Theorem \ref{T:NonProper2}]
Let $\LDS$ be a lifting data set such that the underlying contraction $\om$ is
an isometry and $\om_2\Pi_\cF$ on $\cD_A$ is strongly stable. Then the coefficient
matrix \eqref{CoefMat} for the Redheffer coefficients \eqref{RCLcoefs} associated
with the lifting data set $\LDS$ is unitary according to Theorem \ref{T:RCLRedheffer}.

Moreover, for any $V\in\bS(\cD_A,\cD_{T'})$, since, by Parrott's lemma \cite{P78}, the
operators $B_V$ and $\Ga_V$ in \eqref{RCLsols2}
determine each other uniquely, it follows that $\sV_{B_V}=\sW_{H_V}$,
where $\sV_{B_V}$ and $\sW_{H_V}$ are given by \eqref{nonuniqueSet} and \eqref{nonuniqueSet2},
respectively. Therefore, we obtain Theorem \ref{T:NonProper2} by applying Theorem
\ref{T:NonProper2RED} to the Redheffer representation $\sR_\Phi$.
\end{proof}

\setcounter{equation}{0}
\section{Proof of Theorem \ref{T:MaxPrinciple}}\label{S:maxprinciple}

In this section we derive Theorem \ref{T:MaxPrinciple} from the following
harmonic maximal principle for Redheffer representations.

\begin{theorem}\label{T:MaximumPrinciple}
Let $\Psi_{1,1}$, $\Psi_{1,2}$, $\Psi_{2,1}$ and $\Psi_{2,2}$ be as in
\eqref{RedCoefs2} with $\Psi_{1,1}(0)=0$ and assume that the coefficient
matrix \eqref{RedMat2} is contractive. Then, for the Redheffer
representation $\sR_\Psi$ we have:
\begin{itemize}
\item[(1)] If there exists a $\wtilV\in\bS_0(\cE,\cE')$ with $\|\Ga_{H_{\wtilV}}\|<1$,
then $\|\Ga_{H_{V}}\|<1$ for any $V\in\bS_0(\cE,\cE')$.

\item[(2)] If there exists a $\wtilV\in\bS_0(\cE,\cE')$ with $\|\Ga_{H_{\wtilV}}\|=1$,
then $\|\Ga_{H_{V}}\|=1$ for any $V\in\bS_0(\cE,\cE')$.

\item[(3)] Moreover, for $V\in\bS_0(\cE,\cE')$ and $\wtilV\in\bS(\cE,\cE')$ we have
$\kr D_{\Ga_{H_V}}\subset\kr D_{\Ga_{H_{\wtilV}}}$. In particular,
$\kr D_{\Ga_{H_V}}=\kr D_{\Ga_{H_{\wtilV}}}$
in case $\wtilV,V\in\bS_0(\cE,\cE')$.

\end{itemize}
\end{theorem}

\begin{proof}[\bf Proof]
We first prove (3). Let $V\in\bS_0(\cE,\cE')$ and
$u\in\kr D_{\Ga_{H_V}}$, i.e., $\|\Ga_{H_V}u\|=\|u\|$. Since the operator matrix
$K_V$ of Proposition \ref{P:KVMat2} is contractive, it follows that
$-D_{M_V}(I-M_{\Psi_{1,1}}M_V)^{-1}\Ga_{\Psi_{1,2}}u=0$. The operator
$D_{M_V}(I-M_{\Psi_{1,1}}M_V)^{-1}$ is invertible, because $\|M_V\|<1$.
Thus $\Ga_{\Psi_{1,2}}u=0$. In other words, $\Psi_{1,2}(\la)u=0$ for each
$\la\in\BD$. This implies that $H_{\wtilV}(\la)u=\Psi_{2,2}(\la)u$ for
any $\wtilV\in\bS(\cE,\cE')$ and $\la\in\BD$. Then
$\Ga_{H_{\wtilV}}u=\Ga_{\Psi_{2,2}}u=\Ga_{H_{V}}u$, and thus
$\|\Ga_{H_{\wtilV}}u\|=\|\Ga_{H_{V}}u\|=\|u\|$. Hence $u\in\kr D_{\Ga_{H_{\wtilV}}}$
for any $\wtilV\in\bS(\cE,\cE')$.

Note that (1) is just the contrapositive of (2). To see that (2) holds, we use a similar
argument as in the proof of (3). Assume we have a $V\in\bS_0(\cE,\cE')$ with
$\|\Ga_{H_{V}}\|=1$. Let $u_n$, $n=0,1,\ldots$, be a sequence with $\|u_n\|=1$ and
$\|\Ga_{H_{V}}u_n\|\to1$. Since $K_V$ in Proposition \ref{P:KVMat2} is contractive,
$-D_{M_V}(I-M_{\Psi_{1,1}}M_V)^{-1}\Ga_{\Psi_{1,2}}u_n\to0$, and thus
$\Ga_{\Psi_{1,2}}u_n\to0$; again because $D_{M_V}(I-M_{\Psi_{1,1}}M_V)^{-1}$ is invertible.
Now let $\wtilV$ be an arbitrary function in $\bS_0(\cE,\cE')$.
Using that $\Ga_{H_{\wtilV}}$ and $\Ga_{H_V}$ can be written explicitly as in
\eqref{GaHVexplicit}, it follows that
$\lim_{n\to\infty} \Ga_{H_{\wtilV}}u_n
=\lim_{n\to\infty} \Ga_{\Psi_{2,2}}u_n=\lim_{n\to\infty} \Ga_{H_{V}}u_n$.
In particular, $\|\Ga_{H_{\wtilV}}u_n\|\to1$. Thus $\|\Ga_{H_{\wtilV}}\|=1$.
\end{proof}

\begin{proof}[\bf Proof of Theorem \ref{T:MaxPrinciple}]
Assume there exists a $\wtilV\in\bS_0(\cG,\cD_{\om^*})$ with $\|B_{\wtilV}\|<1$.
Let $H_{\wtilV}=\sR_\Phi[\wtilV]$, with $\Phi$ defined by \eqref{PhiRCL}.
Then both $A$ and $\Ga_{H_{\wtilV}}$ are strict contractions. By Theorem
\ref{T:MaximumPrinciple} it then follows for any $V\in\bS_0(\cG,\cD_{\om^*})$
that the operator $\Ga_{H_{V}}$ is a strict contraction, and thus that $B_V$
is a strict contraction. This proves Part 1. Part 2 follows because it is the
contrapositive of Part 1.

To see that the third statement holds, note that for any $V\in\bS(\cG,\cD_{\om^*})$
the contractive interpolant $B_V$ in \eqref{RCLsols2} satisfies
\[
\|D_{V}h\|=\|D_{\Ga_{H_V}}D_A h\|\quad(h\in\cH).
\]
Thus $h\in\cH$ is in $\kr D_{V}$ if and only if $D_A h$ is in $\kr D_{\Ga_V}$.
The statement then follows by applying Part 3 of Theorem \ref{T:MaximumPrinciple}
to the Redheffer representation $\sR_\Phi$.
\end{proof}

We conclude with the some bounds on $\|\Ga_{H_V}\|$ for $V\in\bS_0(\cE,\cE')$.

\begin{proposition}\label{P:bounds}
Let $\Psi_{1,1}$, $\Psi_{1,2}$, $\Psi_{2,1}$ and $\Psi_{2,2}$ be as in
\eqref{RedCoefs2} with $\Psi_{1,1}(0)=0$ and assume that the coefficient
matrix \eqref{RedMat2} is contractive. For any $V\in\bS_0(\cE,\cE')$ and $u\in\cU$ we have
\begin{equation}\label{bound1}
\|\Ga_{H_V}u\|^2\leq\|u\|^2-\frac{1-\|M_V\|^2}{\|I-M_{\Psi_{1,1}}M_V\|^2}\|\Ga_{\Psi_{1,2}}u\|^2.
\end{equation}
In particular,
\begin{equation}\label{bound2}
\|\Ga_{H_V}u\|^2\leq\|u\|^2-\frac{1-\|M_V\|}{1+\|M_V\|}\|\Ga_{\Psi_{1,2}}u\|^2.
\end{equation}
If, in addition, the coefficient matrix \eqref{RedMat2} is an isometry, then
\begin{equation}\label{bound3}
\|\Ga_{H_V}u\|^2\leq\frac{2\|M_V\|}{1+\|M_V\|}\|u\|^2
+\frac{1-\|M_V\|}{1+\|M_V\|}\|\Ga_{\Psi_{2,2}}u\|^2.
\end{equation}
\end{proposition}

\begin{proof}[\bf Proof]
{}From the fact that the operator matrix $K_V$ in \eqref{KV2} is contractive we obtain
\[
\|\Ga_{H_V}u\|^2\leq \|u\|^2-\|D_V(I-M_{\Psi_{11}}M_V)^{-1}\Ga_{\Psi_{1,2}}u\|^2.
\]
Using that for an invertible operator $T\in\cL(\cH)$ and $h\in\cH$
we have the inequality $\|Th\|\geq\|T^{-1}\|^{-1}\|h\|$, it follows that
\[
\|\Ga_{H_V}u\|^2\leq
\|u\|^2-\|D_V^{-1}\|^{-2}\|(I-M_{\Psi_{11}}M_V)\|^{-2}\|\Ga_{\Psi_{1,2}}u\|^2.
\]
The first inequality for $\|\Ga_{H_V}u\|$ holds because $\|D_{M_V}^{-1}\|\leq(1-\|M_V\|^2)^{-\half}$.
Indeed, for any $h\in H^2(\cG)$ we have $\|D_{M_V}h\|^2=\|h\|^2-\|M_Vh\|^2\geq(1-\|M_V\|^2)\|h\|^2$.
So for $h=D_{M_{V}}^{-1}k$ we obtain $\|k\|^2\geq (1-\|M_V\|^2)\|D_{M_{V}}^{-1}k\|^2$, which proves
our claim.

To see that the second inequality holds, we use
\[
\|I-M_{\Psi_{1,1}}M_V\|\leq 1+\|M_{\Psi_{1,1}}M_V\|\leq1+\|M_V\|,
\]
in order to obtain
\[
\frac{1-\|M_V\|^2}{\|I-M_{\Psi_{1,1}}M_V\|^2}\geq
\frac{1-\|M_V\|^2}{(1+\|M_V\|)^2}=\frac{1-\|M_V\|}{1+\|M_V\|}.
\]

In case the coefficient matrix $K_0$ is an isometry,
$\|\Ga_{\Psi_{1,2}}u\|^2=\|u\|^2-\|\Ga_{\Psi_{2,2}}u\|^2$. Hence,
starting from the second inequality,
\begin{eqnarray*}
\|\Ga_{H_V}u\|^2&\leq&
\|u\|^2-\frac{1-\|M_V\|}{1+\|M_V\|}(\|u\|^2-\|\Ga_{\Psi_{2,2}}u\|^2)\\
&=&\frac{2\|M_V\|}{1+\|M_V\|}\|u\|^2+\frac{1-\|M_V\|}{1+\|M_V\|}\|\Ga_{\Psi_{2,2}}u\|^2.
\end{eqnarray*}
\end{proof}


\end{document}